\documentclass{article}%
\usepackage{amsfonts}
\usepackage{graphicx}
\usepackage{amsmath}
\usepackage{amssymb}%
\providecommand{\U}[1]{\protect\rule{.1in}{.1in}}
\newtheorem {theorem}{Theorem}[section]
\newtheorem {proposition}{Proposition}[section]
\newtheorem {corollary}{Corollary}[section]

\newtheorem{definition}{Definition}[section]

\newtheorem{lemma}{Lemma}[section]

\newtheorem{remark}{Remark}[section]

\begin{document}

\begin{center}
{\LARGE Central Limit Theorem for Linear Processes}

\bigskip

{\LARGE with Infinite Variance }\vskip15pt

\bigskip Magda Peligrad$^{a}$\footnote{Supported in part by a Charles Phelps
Taft Memorial Fund grant and NSA\ grant H98230-09-1-0005} and Hailin
Sang$^{b}$

$^{a}$ Department of Mathematical Sciences, University of Cincinnati, PO Box
210025, Cincinnati, OH 45221-0025, USA. E-mail address: peligrm@ucmail.uc.edu

\bigskip

$^{b}$ National Institute of Statistical Sciences, PO Box 14006, Research
Triangle Park, NC 27709, USA. E-mail address: sang@niss.org

\bigskip
\end{center}

\textit{MSC 2000 subject classification}: Primary: 60F05, 60G10, 60G42

Key words and phrases: linear process, central limit theorem, martingale,
mixing, infinite variance.

\bigskip

\textbf{Abstract}

\bigskip

This paper addresses the following classical question: giving a sequence of
identically distributed random variables in the domain of attraction of a
normal law, does the associated linear process satisfy the central limit
theorem? We study the question for several classes of dependent random
variables. For independent and identically distributed random variables we
show that the central limit theorem for the linear process is equivalent to
the fact that the variables are in the domain of attraction of a normal law,
answering in this way an open problem in the literature. The study is also
motivated by models arising in economic applications where often the
innovations have infinite variance, coefficients are not absolutely summable,
and the innovations are dependent.

\section{Introduction and notations}

Let $(\xi_{n})_{n\in\mathbb{Z}}$ be a sequence of identically distributed
random variables and let $(c_{ni})_{1\leq i\leq m_{n}}$ be a triangular array
of numbers. In this paper we analyze the asymptotic behavior of statistics of
the type
\begin{equation}
S_{n}=\sum_{i=1}^{m_{n}}c_{ni}\xi_{i} \label{gln}%
\end{equation}
when the variables are centered and satisfy:%
\begin{equation}
H(x)=\mathbb{E(}\xi_{0}^{2}I(|\xi_{0}|\leq x))\text{ is a slowly varying
function at }\infty\text{.} \label{DA}%
\end{equation}
This tail condition is highly relevant to the central limit theory. For
independent, identically distributed, centered variables this condition is
equivalent to the fact that the variables are in the domain of attraction of
the normal law. This means: there is a sequence of constants $b_{n}%
\rightarrow\infty$ such that $\sum_{i=1}^{n}\xi_{i}/b_{n}$ is convergent in
distribution to a standard normal variable (see for instance Feller, 1966;
Ibragimov and Linnik, 1971; Araujo and Gin\'{e}, 1980). It is an open question
to extend the general central limit theorem from equal weights to weighted
sums of i.i.d. random variables with infinite variance. Linear combinations of
identically distributed random variables are important, since many random
evolutions and also statistical procedures such as parametric or nonparametric
estimation of regression with fixed design, produce statistics of type
(\ref{gln}) (see for instance Chapter 9 in Beran, 1994, for the case of
parametric regression, or the paper by Robinson, 1997, where kernel estimators
are used for nonparametric regression). One example is the simple parametric
regression model $Y_{i}=\beta\alpha_{i}+\xi_{i}$ where $(\xi_{i})$ is a
sequence of identically distributed random variables with marginal
distribution satisfying (\ref{DA}), $(\alpha_{i})$ is a sequence of real
numbers and \ $\beta$ is the parameter of interest. The least squares
estimator $\hat{\beta}$ of $\ \beta$, based on a sample of size $n,$ satisfies
$S_{n}=\hat{\beta}-\ \beta=(\sum_{i=1}^{n}\alpha_{i}^{2})^{-1}\sum_{i=1}%
^{n}\alpha_{i}\xi_{i}$. So the representation of type (\ref{gln}) holds with
$c_{ni}=\alpha_{i}/(\sum_{i=1}^{n}\alpha_{i}^{2})$.

We shall also see that the asymptotic behavior of the sums of variables of the
form
\begin{equation}
X_{k}=\sum_{j=-\infty}^{\infty}a_{k+j}\xi_{j} \label{ln}%
\end{equation}
can be obtained by studying sums of the type (\ref{gln}). We shall refer to
such a process as to a linear process with innovations $(\xi_{i}%
)_{i\in{\mathbb{Z}}}.$ In 1971 Ibragimov and Linnik extended the central limit
theorem from i.i.d. random variables to linear processes defined by (\ref{ln})
for innovations that have finite second moment, under the conditions
$\sum_{j=-\infty}^{\infty}a_{k}^{2}<\infty$ and $stdev(\sum_{j=1}^{n}%
X_{j})\rightarrow\infty.$ They showed that
\[
\sum_{j=1}^{n}X_{j}/stdev(\sum_{j=1}^{n}X_{j})\overset{D}{\rightarrow}N(0,1).
\]
This result is striking, since $var(\sum_{j=1}^{n}X_{j})$ can be of order
different of $n$; practically it can be any positive sequence going to
infinite of an order $o(n^{2}).$ It was conjectured that a similar result
might hold without the assumption of finite second moment. Steps in this
direction are papers by Knight (1991), Mikosch et al. (1995) and Wu (2003) who
studied this problem under the additional assumption $\sum_{j=-\infty}%
^{\infty}|a_{k}|<\infty.$ Our Theorem \ref{Equivalence} positively answers
this conjecture. Under condition (\ref{DA}) we show that $X_{k}$ is well
defined, if and only if%
\[
\sum_{j\in\mathbb{Z},a_{j}\neq0}a_{j}^{2}H(|a_{j}|^{-1})<\infty,
\]
and under this condition we show that the central limit theorem for
$\sum_{j=1}^{n}X_{j}$ properly normalized is equivalent to condition (\ref{DA}).

As an example in this class we mention the particular linear process with
regularly varying weights with exponent $\alpha$ where $1/2<\alpha<1.$ This
means that the coefficients are of the form $a_{n}=n^{-\alpha}L(n)$, where for
$n\geq1$, $a_{n}=0$ for $n\leq0,$ and $L(n)$ is a slowly varying function at
$\infty$. It incorporates the fractionally integrated processes that play an
important role in financial econometrics, climatology and so on and they are
widely studied. Such processes are defined for $0<d<1/2$ by
\[
X_{k}=(1-B)^{-d}\xi_{k}=\sum_{i\geq0}a_{i}\xi_{k-i}\text{ where }a_{i}%
=\frac{\Gamma(i+d)}{\Gamma(d)\Gamma(i+1)}%
\]
and $B$ is the backward shift operator, $B\varepsilon_{k}=\varepsilon_{k-1}$.
For this example, by the well known fact that for any real $x,$ $\lim
_{n\rightarrow\infty}\Gamma(n+x)/n^{x}\Gamma(n)=1$ we have$\ \lim
_{n\rightarrow\infty}a_{n}/n^{d-1}=1/\Gamma(d)$. Notice that these processes
have long memory because $\sum_{j\geq1}|a_{j}|=\infty.$ This particular class
was recently investigated by Peligrad and Sang (2010), where further reaching
properties were pointed out.

Our study is not restricted to the class of independent identically
distributed random variables. We also consider larger classes including
martingales and mixing processes. The results obtained for the class of
martingale innovations are also useful to study more general innovations that
can be approximated by martingale differences. The martingale approximation
method was recently used by Jara et al (2009) to study the attraction to
stable laws with exponent $\alpha,$ $\alpha\in(0,2)$ for additive functionals
of a stationary Markov chain.

There is a huge literature on the central limit theorem for linear processes
with dependent innovations and finite second moment but we are not aware of
any study considering the infinite variance case in its full generality. A
step in this direction, under the assumption $\sum_{j=-\infty}^{\infty}%
|a_{k}|<\infty,$ is the paper by Tyran-Kami\'{n}ska (2010).

In all the central limit theorems for variables with infinite variance the
construction of the normalizer is rather complicated and is based heavily on
the function $H(x).$ This is the reason why it is important to replace the
normalizer by a selfnormalizer, constructed from the data. We mention in this
direction the recent results by Kulik (2006), under the assumption that
$\sum_{j=-\infty}^{\infty}|a_{k}|<\infty$ and by Peligrad and Sang (2010) for
regularly varying weights with exponent $\alpha$ where $1/2<\alpha<1$. In this
paper, as in Mason (2005), we suggest a Raikov type selfnormalizer based on a
weighted sum of squares of the innovations.

Our paper is organized in the following way: Section 2 contains the
definitions and the results, Section 3 contains the proofs. For convenience,
in the Appendix, we give some auxiliary results and we also mention some known
facts needed for the proofs.

In the sequel we shall use the following notations: a double indexed sequence
with indexes $n$ and $i$ will be denoted by $a_{ni}$ and sometimes $a_{n,i};$
we use the notation $a_{n}\sim b_{n}$ instead of $a_{n}/b_{n}\rightarrow1$;
$a_{n}=o(b_{n})$ means that $a_{n}/b_{n}\rightarrow0;$ $I(A)$ denotes the
indicator function of $A$; the notation $\Rightarrow$ is used for convergence
in distribution and also for convergence in probability to a constant.

In this paper we shall make two conventions in order to simplify the notations.

\textbf{Convention 1}. By convention, for $x=0$, $|x|H(|x|^{-1})=0.$ For
instance we can write instead $\sum_{j\in\mathbb{Z},\text{ }a_{j}\neq0}%
a_{j}^{2}H(|a_{j}|^{-1})<\infty,$ simply $\sum_{j\in\mathbb{Z}}a_{j}%
^{2}H(|a_{j}|^{-1})<\infty.$

\textbf{Convention 2.} The second convention refers to the function $H(x)$
defined in (\ref{DA})$.$ Since the case $\mathbb{E(}\xi_{0}^{2})<\infty$ is
known, we shall consider the case $\mathbb{E(}\xi_{0}^{2})=\infty.$ Let
$b=\inf\left\{  x\geq0:H(x)>1\right\}  \ $\ and $H_{b}(x)=H(x\vee(b+1)).$ Then
clearly $b<\infty$, $H_{b}(x)\geq1$ and $H_{b}(x)=H(x)$ for $x>b+1.$ From now
on we shall redenote $H_{b}(x)$ by $H(x).$ Therefore, since our results are
asymptotic, without restricting the generality we shall assume that
$H(x)\geq1$ for all $x\geq0.$

\section{Results}

Our first results treat the general weights and identically distributed
martingale differences with infinite second moment. The case of finite second
moment was treated in Peligrad and Utev (1997, 2006).

We shall establish first a general theorem for martingale differences under a
convergence in probability condition (\ref{convprob}). This condition will be
verified in the next main results for classes of martingale differences and
i.i.d. random variables.

\begin{theorem}
\label{genthm} Let $(\xi_{k})_{k\in\mathbb{Z}}$ be a sequence of identically
distributed martingale differences adapted to the filtration $(\mathcal{F}%
_{k})_{k\in\mathbb{Z}}$ that satisfy (\ref{DA}) and let $(c_{nk})_{1\leq k\leq
m_{n}}$ be a triangular array of real numbers, such that%
\begin{equation}
\text{ }\sup_{n}\sum_{k=1}^{m_{n}}c_{nk}^{2}H(|c_{nk}|^{-1})<\infty\text{ and
}\max_{1\leq k\leq m_{n}}c_{nk}^{2}H(|c_{nk}|^{-1})\rightarrow0\;\text{as}%
\;n\rightarrow\infty. \label{gen}%
\end{equation}
Assume%
\begin{equation}
\sum_{k=1}^{m_{n}}c_{nk}^{2}\xi_{k}^{2}\Rightarrow1\text{ as}\;\ n\rightarrow
\infty\text{.} \label{convprob}%
\end{equation}
Then
\begin{equation}
\sum_{k=1}^{m_{n}}c_{nk}\xi_{k}\Rightarrow N(0,1)\;\;\text{as}\;\;n\rightarrow
\infty\text{.} \label{CLT}%
\end{equation}

\end{theorem}

We shall mention two pairwise mixing type conditions that are sufficient for
(\ref{convprob}).

\begin{proposition}
\label{mix} Assume that all the conditions of Theorem \ref{genthm} (except for
(\ref{convprob})) are satisfied. Assume that one of the following two
conditions holds:\newline($M_{1}$) There is a sequence of positive numbers
$\psi_{k}\rightarrow0$ such that for all $a$ and $b$ positive numbers and all
integers $j,$
\[
cov(\xi_{j}^{2}I(|\xi_{0}|\leq a),\xi_{j+k}^{2}I(|\xi_{k}|\leq b))\leq\psi
_{k}\mathbb{E(}\xi_{0}^{2}I(|\xi_{0}|\leq b))\mathbb{E(}\xi_{0}^{2}I(|\xi
_{0}|\leq a))\text{.}%
\]
($M_{2}$) There is a sequence of positive numbers $(\varphi_{k})\ $with
$\sum_{k\geq1}^{\ }\varphi_{k}<\infty$, such that for all $a$ and $b$ and all
integers $j,$ \
\[
cov(\xi_{j}^{2}I(|\xi_{0}|\leq a),\xi_{j+k}^{2}I(|\xi_{k}|\leq b))\leq
a^{2}\varphi_{k}\mathbb{E(}\xi_{0}^{2}I(|\xi_{0}|\leq b))\text{.}%
\]
If either ($M_{1}$) or ($M_{2}$) holds then (\ref{convprob}) is satisfied and
therefore the conclusion of Theorem \ref{genthm} holds.
\end{proposition}

\begin{remark}
According to the above proposition we mention that for independent identically
distributed innovations satisfying (\ref{DA}) and \ coefficients $(c_{ni})$
satisfying condition (\ref{gen}) the central limit theorem (\ref{CLT}) holds.
\end{remark}

For further applications to time series we shall comment on the normalized
form of the above results, which is important for the case when condition
(\ref{gen}) is not satisfied.

Recall Conventions 1 and 2 and define:
\begin{equation}
D_{n}=\inf\left\{  s\geq1:\sum_{k=1}^{m_{n}}\frac{c_{nk}^{2}}{s^{2}}H\left(
\frac{s}{c_{nk}}\right)  \leq1\right\}  \text{.} \label{defB}%
\end{equation}
$D_{n}$ is well defined since by Convention 2, $H(x)\geq1$ for all $x$ and,
since $H(x)$ is a slowly varying function, we have $\lim_{x\rightarrow\infty
}x^{-2}H(x)=0$. By using this definition along with Theorem \ref{genthm} we
obtain the following corollary:

\begin{corollary}
\label{CLTnorm}Let $(\xi_{k})_{k\in\mathbb{Z}}$ \ be as in Theorem
\ref{genthm} and assume
\begin{equation}
\max_{1\leq k\leq m_{n}}\frac{c_{nk}^{2}}{D_{n}^{2}}H\left(  \frac{D_{n}%
}{|c_{nk}|}\right)  \rightarrow0\;\text{as}\;n\rightarrow\infty\label{coeffD}%
\end{equation}
and
\begin{equation}
\frac{1}{D_{n}^{2}}\sum_{k=1}^{m_{n}}c_{nk}^{2}\xi_{k}^{2}\Rightarrow1\text{
as}\;\;n\rightarrow\infty\text{.} \label{condconv}%
\end{equation}
Then
\begin{equation}
\frac{1}{D_{n}}\sum_{k=1}^{m_{n}}c_{nk}\xi_{k}\Rightarrow N(0,1)\;\;\text{as}%
\;\;n\rightarrow\infty\text{.} \label{CLTNor}%
\end{equation}
Moreover, as in Proposition \ref{mix}, condition (\ref{condconv}) is satisfied
under either (M$_{1}$) or (M$_{2}$).
\end{corollary}

Clearly, by combining (\ref{condconv}) with (\ref{CLTNor}) we also have%

\begin{equation}
\frac{1}{(\sum_{k=1}^{m_{n}}c_{nk}^{2}\xi_{k}^{2})^{1/2}}\sum_{k=1}^{m_{n}%
}c_{nk}\xi_{k}\Rightarrow N(0,1)\;\;\text{as}\;\;n\rightarrow\infty\text{.}
\label{selfnorm}%
\end{equation}
For equal weights, when $c_{n,i}=1$ for all $n$ and $1\leq i\leq n,$ $D_{n}$
becomes the standard normalizer for the central limit theorem for variables in
the domain of attraction of a normal law:
\[
D_{n}=\inf\left\{  s\geq1:\frac{n}{s^{2}}H\left(  s\right)  \leq1\right\}
\text{.}%
\]
Then, it is well known that $D_{n}\rightarrow\infty$ and, by the properties of
slowly varying function, condition (\ref{coeffD}) is satisfied with $m_{n}=n.$
For this case we easily obtain the following central limit theorem for
martingales with infinite variance:

\begin{corollary}
Let $(\xi_{k})_{k\in\mathbb{Z}}$ be a sequence of identically distributed
martingale differences adapted to the filtration $(\mathcal{F}_{k}%
)_{k\in\mathbb{Z}}$ satisfying (\ref{DA}) and assume
\begin{equation}
\frac{1}{D_{n}^{2}}\sum_{k=1}^{n}\xi_{k}^{2}\Rightarrow1\text{ \ as}%
\;\;n\rightarrow\infty\text{.} \label{condconv1}%
\end{equation}
Let $M_{n}=\sum_{k=1}^{n}\xi_{k}.$ Then
\[
\frac{M_{n}}{D_{n}}\Rightarrow N(0,1)\;\;\text{as}\;\;n\rightarrow
\infty\text{.}%
\]
Moreover, as in Proposition \ref{mix}, condition (\ref{condconv1}) is
satisfied under either (M$_{1}$) or (M$_{2}$).
\end{corollary}

\begin{remark}
\label{Gordin}This corollary can be used to obtain the CLT for classes of
stochastic processes that can be approximated by stationary martingale
differences. For instance, assume $(\eta_{k})_{k\in\mathbb{Z}}$ is a
stationary sequence of centered, integrable random variables, $\mathcal{F}%
_{i}=\sigma(\eta_{j},j\leq i)$ and $V_{n}=\sum_{k=1}^{n}\eta_{k}.$ An idea,
going back to Gordin (1969), is to decompose $V_{n}$ into a martingale with
stationary differences and a telescoping rest called coboundary. More
precisely, $\eta_{n}=\xi_{n}+Z_{n-1}-Z_{n},$ where $Z_{n}$ is a stationary
integrable sequence and $\xi_{n}$ is a stationary sequence of martingale
differences. Voln\'{y} (1993) gave necessary and sufficient conditions for
such an approximation. Under the assumption $\mathbb{E}(V_{n}|\mathcal{F}%
_{0})$ is convergent in $\mathbb{L}_{1},$ we have $V_{n}=M_{n}+R_{n}$ where
$M_{n}$ is a martingale with stationary differences adapted to $\mathcal{F}%
_{n}$, and $R_{n}=Z_{0}-Z_{n}$ has the property that $E|R_{n}|<\infty.$ Then
clearly, after normalizing by a sequence of constants $B_{n}$ converging to
$\infty$, the limiting distribution of $V_{n}/B_{n}$ is equivalent to the
limiting distribution of $M_{n}/B_{n}.$
\end{remark}

By combining our Corollary \ref{CLTnorm} with the main result in Gin\'{e} et
al. (1997) we formulate next a theorem for i.i.d. sequences.

\begin{theorem}
Let $(\xi_{k})_{k\in\mathbb{Z}}$ be a sequence of independent and identically
distributed centered random variables. Then the following three statements are
equivalent:\newline(1) $\xi_{0}$ is in the domain of attraction of a normal
law ( i.e. condition (\ref{DA}) is satisfied).\newline(2) For any sequence of
constants $(c_{nk})_{1\leq k\leq m_{n}}$ satisfying (\ref{coeffD}) the CLT in
( \ref{CLTNor}) holds.\newline(3) For any sequence of constants $(c_{nk}%
)_{1\leq k\leq m_{n}}$ satisfying (\ref{coeffD}) the selfnormalized CLT in
(\ref{selfnorm}) holds.
\end{theorem}

The implication (2)$\rightarrow$(3) also follows by Mason (2005). We shall
apply our general results for time series of the form (\ref{ln}) with
identically distributed martingale differences innovations. We shall prove
first the following proposition:

\begin{proposition}
\label{Exist}Let $(\xi_{k})_{k\in\mathbb{Z}}$ be a sequence of identically
distributed martingale differences adapted to the filtration $(\mathcal{F}%
_{k})_{k\in\mathbb{Z}}$ that satisfy (\ref{DA}). The linear process
$X_{0}={\sum_{j=-\infty}^{\infty}}a_{j}\xi_{j}$ is well defined in the almost
sure sense under the condition
\begin{equation}
\sum_{j\in\mathbb{Z}}a_{j}^{2}H(|a_{j}|^{-1})<\infty. \label{coeff0}%
\end{equation}
If the innovations are independent and identically distributed random
variables satisfying (\ref{DA}) then condition (\ref{coeff0}) is necessary and
sufficient for the existence of $X_{0}$ a.s$.$
\end{proposition}

Denote
\[
b_{nj}=a_{j+1}+\cdots+a_{j+n}%
\]
and with this notation $\ $%
\begin{equation}
S_{n}=\sum_{k=1}^{n}X_{k}=\sum_{j\in{\mathbb{Z}}}b_{nj}\xi_{j}\text{.}
\label{defS}%
\end{equation}
Construct $D_{n}$ by (\ref{defB}) where we replace $c_{nj}$ by $b_{nj}$. Then
we have:

\begin{theorem}
\label{linCLT}Assume that $(\xi_{k})_{k\in\mathbb{Z}}$ is a sequence of
identically distributed martingale differences satisfying (\ref{DA}), the
coefficients satisfy condition (\ref{coeff0}), and $\sum_{k}b_{nk}%
^{2}\rightarrow\infty$. Assume in addition that
\begin{equation}
\frac{1}{D_{n}^{2}}\sum_{k}b_{nk}^{2}\xi_{k}^{2}\Rightarrow1\text{
as}\;\;n\rightarrow\infty\text{.} \label{condconv2}%
\end{equation}
Then%
\begin{equation}
\frac{S_{n}}{D_{n}}\Rightarrow N(0,1)\;\;\text{as}\;\;n\rightarrow
\infty\text{.} \label{CLT2}%
\end{equation}

\end{theorem}

Notice that the normalizer $D_{n}$ is rather complicated and contains the
slowly varying function $H(x).$ By combining however the convergences in
(\ref{condconv2}) and (\ref{CLT2}) it is easy to see that we can use in
applications the selfnormalized form of this theorem.
\begin{equation}
\frac{S_{n}}{(\sum_{k}b_{nk}^{2}\xi_{k}^{2})^{1/2}}\Rightarrow
N(0,1)\;\;\text{as}\;\;n\rightarrow\infty\text{.} \label{self}%
\end{equation}
By simple arguments $D_{n}$ in (\ref{CLT2}) can also be replaced by $\sqrt
{\pi/2}\mathbb{E(}|S_{n}|),$ or a consistent estimator of this quantity.

We mention now sufficient conditions for the validity of (\ref{condconv2}).

\begin{theorem}
\label{linmart}Assume that $(\xi_{k})_{k\in\mathbb{Z}}$ is a sequence of
identically distributed martingale differences that satisfy (\ref{DA}), the
coefficients $(a_{i})$ are as in Theorem \ref{linCLT} and one of the
conditions (M$_{1}$) or (M$_{2}$) of Proposition \ref{mix} is satisfied. Then
both the central limit theorem (\ref{CLT2}) and its selfnormalized form
(\ref{self}) hold.
\end{theorem}

In the independent case, by combining theorem (\ref{linmart}) with the result
on selfnormalized CLT\ in Gin\'{e} et al. (1997), we have:

\begin{theorem}
\label{Equivalence}Let $(\xi_{k})_{k\in\mathbb{Z}}$ be a sequence of
independent and identically distributed centered random variables. Then the
following three statements are equivalent:\newline(1) $\xi_{0}$ satisfies
condition (\ref{DA}).\newline(2) For any sequence of constants $(a_{n}%
)_{n\in\mathbb{Z}}$ satisfying (\ref{coeff0}) and $\sum_{k}b_{nk}%
^{2}\rightarrow\infty$ the CLT in (\ref{CLT2}) holds.\newline(3) For any
sequence of constants $(a_{n})_{n\in\mathbb{Z}}$ satisfying (\ref{coeff0}) and
$\sum_{k}b_{nk}^{2}\rightarrow\infty$ the selfnormalized CLT in (\ref{self}) holds.
\end{theorem}

This theorem is an extension of Theorem 18.6.5 in Ibragimov and Linnik (1971)
from i.i.d. innovations with finite second moment to innovations in the domain
of attraction of a normal law. It positively answers the question on the
stability of the central limit theorem for i.i.d. variables under formation of
linear sums. The implication (2)$\rightarrow$(3) also follows by Mason (2005).

The casual linear process is obtained when $a_{i}=0$ for $i\leq0.$ Then
$X_{0}={\sum_{j=-\infty}^{0}}a_{j}\xi_{j}$ is well defined if and only if
\begin{equation}
\sum_{j\geq0}a_{j}^{2}H(|a_{j}|^{-1})<\infty. \label{coeff02}%
\end{equation}
For this case the coefficients $b_{nj}$ have the following expression
\begin{align*}
b_{nj}  &  =a_{1}+\ldots+a_{j}\text{ for }j<n\\
b_{nj}  &  =a_{j-n+1}+\ldots+a_{j}\,\text{\ for }j\geq n
\end{align*}
and with this notation,%
\begin{equation}
S_{n}=\sum_{i=1}^{\infty}b_{ni}\mathbb{\xi}_{n-i}.\nonumber
\end{equation}
For particular casual linear processes with coefficients $a_{i}=i^{-\alpha
}L(i)$, where $1/2<\alpha<1$ and $L(i)$ is a slowly varying function at
$\infty$ in the strong sense (i.e. there is $h(t)$ continuous such that
$L(n)=h(n)$ and $h(t)$ is slowly varying), the normalizer can be made more
precise. Peligrad and Sang (2010) studied the case of i.i.d. innovations and
showed that%
\[
D_{n}^{2}\sim c_{\alpha}H(\eta_{n})n^{3-2\alpha}L^{2}(n)
\]
where $c_{\alpha}=(1-\alpha)^{-2}\int_{0}^{\infty}[x^{1-\alpha}-\max
(x-1,0)^{1-\alpha}]^{2}dx$ and $\eta_{n}\ $is defined by
\begin{equation}
\eta_{j}=\inf\left\{  s>1:\frac{H(s)}{s^{2}}\leq\frac{1}{j}\right\}
,\;\;\;j=1,2,\cdots\label{defeta}%
\end{equation}
Furthermore, in this context, Peligrad and Sang (2010) showed that the
selfnormalizer can be estimated by observing only the variables $X_{k}$. More
precisely they showed that $c_{\alpha}n^{2}a_{n}^{2}V_{n}^{2}/(A^{2}D_{n}%
^{2})\Rightarrow1$ , where $V_{n}^{2}=\sum_{i=1}^{n}X_{i}^{2}$ and $A^{2}%
=\sum_{i=1}^{\infty}a_{i}^{2}.$ Therefore
\[
\frac{S_{n}}{na_{n}V_{n}}\Rightarrow N(0,\frac{c_{\alpha}}{A^{2}}).
\]
By combining this result with Theorem \ref{linCLT} we notice that for this
case we have the following striking law of large numbers:
\[
\frac{A^{2}\sum_{i=1}^{\infty}b_{ni}^{2}\mathbb{\xi}_{n-i}^{2}}{c_{\alpha
}n^{2}a_{n}^{2}V_{n}^{2}}\Rightarrow1\text{ as }n\rightarrow\infty.
\]
Another particular case of linear processes with i.i.d. innovations in the
domain of attraction of a normal law was studied by Kulik (2006), under the
condition $\sum_{j>0}|a_{j}|<\infty.$ For this case his result is:%
\begin{equation}
\frac{S_{n}}{V_{n}}\Rightarrow N(0,\frac{|\sum_{i>0}a_{i}|^{2}}{A^{2}}).
\label{short}%
\end{equation}
It is an interesting related question to extend Kulik's result to martingale
differences. We shall not pursue this path here where the goal is to consider
general coefficients.

We shall also study the case of weak dependent random variables whose
definition is based on the maximum coefficient of correlation.

\begin{definition}
Let $\mathcal{A}$ and $\mathcal{B}$ be two $\sigma$-algebras of events and
define
\[
\rho(\mathcal{A}\text{,}\mathcal{B})=\sup_{f\in\mathbb{L}_{2}(\mathcal{A}%
),g\in\mathbb{L}_{2}(\mathcal{B})}|corr(f,g)|
\]
where $\mathbb{L}_{2}(\mathcal{A})$ denotes the the class of random variables
that are $\mathcal{A}$-measurable and square integrable.
\end{definition}

\begin{definition}
Let $(\xi_{k})_{k\in\mathbb{Z}}$ be a sequence of random variables and let
$\mathcal{F}_{n}^{m}=\sigma(\xi_{i},n\leq i\leq m)$. We call the sequence
$\rho$-mixing if%
\[
\rho_{n}=\sup_{k}\rho(\mathcal{F}_{1}^{k}\text{,}\mathcal{F}_{k+n}^{\infty
})\rightarrow0\text{ as }n\rightarrow\infty.
\]

\end{definition}

This class is significant for studying functions of Gaussian processes, as
well as additive functionals of Markov processes. A convenient reference for
basic properties and the computation of these coefficients for functions of
Markov chains and functions of Gaussian processes is Bradley (2007), chapters
7, 9 and 27.

The next theorem solves the same problem as Theorem \ref{genthm} for this
class of dependent random variables. The conditions imposed to the variables
and mixing rates are similar to those used by Bradley (1988) who studied the
central limit theorem for partial sums of stationary $\rho$-mixing sequences
under (\ref{DA}). Bradley's result was extended in Shao (1993) in several
directions, but still for partial sums. Our theorem extends Theorem 1 of
Bradley (1988) from equal weights to linear processes and Theorem 2.2 (b) in
Peligrad and Utev (1997) to variables with infinite second moment.

\begin{theorem}
\label{MIXING}Let $(\xi_{k})$ be a sequence of centered identically
distributed random variables satisfying (\ref{DA}). Assume that $(\xi_{k})$ is
$\rho$-mixing with $\sum_{k}\rho(2^{k})<\infty$ and $\rho(1)<1$.\newline Let
$(c_{nk})$ be a triangular array of real numbers satisfying
\[
\sup_{n}\sum_{k=1}^{m_{n}}c_{nk}^{2}H(|c_{nk}|^{-1})<\infty\text{ and }%
\max_{1\leq k\leq m_{n}}c_{nk}^{2}H(|c_{nk}|^{-1})\rightarrow0\;\text{as}%
\;n\rightarrow\infty.
\]
Then
\[
\frac{1}{B_{n}}\sum_{k=1}^{m_{n}}c_{nk}\xi_{k}\Rightarrow N(0,1)\;\;\text{as}%
\;\;n\rightarrow\infty\text{,}%
\]
where $B_{n\text{ }}=(\pi/2)^{1/2}$ $\mathbb{E}|\sum_{k=1}^{m_{n}}c_{nk}%
\xi_{k}|$ .
\end{theorem}

\section{Proofs}

{\Large Proof of Theorem \ref{genthm}.}

\bigskip

The proof of Theorem \ref{genthm} involves a few steps. Define%
\[
\xi_{i}^{\prime}=\xi_{ni}^{\prime}=\xi_{i}I(|c_{ni}\xi_{i}|\leq1)-\mathbb{E}%
_{i-1}(\xi_{i}I(|c_{ni}\xi_{i}|\leq1))
\]
and%
\[
\xi_{i}^{\prime\prime}=\xi_{ni}^{\prime\prime}=\xi_{i}I(|c_{ni}\xi
_{i}|>1)-\mathbb{E}_{i-1}(\xi_{i}I(|c_{ni}\xi_{i}|>1))\text{,}%
\]
where we used the notation $\mathbb{E}_{i}(X)$ instead of $\mathbb{E}%
(X|\mathcal{F}_{i})$.

We show now that
\[
\sum_{k=1}^{m_{n}}c_{nk}\xi_{k}^{\prime\prime}\Rightarrow0\text{.}%
\]
This is true because by the item 2 of Lemma \ref{csorgo} and (\ref{gen}),
\begin{gather}
\mathbb{E}|\sum_{k=1}^{m_{n}}c_{nk}\xi_{k}^{\prime\prime}|\leq2\sum
_{i=1}^{m_{n}}|c_{ni}|\mathbb{E(}|\xi_{0}|I(|c_{ni}\xi_{0}|>1))\label{negl}\\
=o(\sum_{i=1}^{m_{n}}c_{ni}^{2}H(|c_{ni}|^{-1})=o(1)\;\;\text{as}%
\;\;n\rightarrow\infty\text{.}\nonumber
\end{gather}
To prove the theorem, by Theorem 3.1 in Billingsley (1999) and (\ref{negl}),
it is enough to study the limiting distribution for the linear process
associated to $(\xi_{i}^{\prime}).$

We shall verify the sufficient conditions for the CLT for sums of a triangular
array of martingale differences with finite second moment, given for
convenience in Theorem \ref{MARTCLT} in the Appendix. We start by verifying
the point $(a)$ of Theorem \ref{MARTCLT}. Fix $0<\varepsilon<1$ and notice
that by the properties of conditional expectations and item 4 in Lemma
\ref{csorgo}, we have:%
\begin{gather*}
\mathbb{E}(\max_{1\leq i\leq m_{n}}|c_{ni}\xi_{i}^{\prime}|^{2})\leq
\varepsilon^{2}+\sum_{k=1}^{m_{n}}c_{nk}^{2}\mathbb{E(}(\xi_{k}^{\prime}%
)^{2}I(|c_{nk}\xi_{k}^{\prime}|>\varepsilon))\\
\leq\varepsilon^{2}+\frac{1}{\varepsilon}\sum_{k=1}^{m_{n}}|c_{nk}%
|^{3}\mathbb{E}|\xi_{k}^{\prime}|^{3}\leq\varepsilon^{2}+\frac{8}{\varepsilon
}\sum_{k=1}^{m_{n}}|c_{nk}|^{3}\mathbb{E(}|\xi_{k}|^{3}I(|c_{nk}\xi_{k}%
|\leq1))\\
\leq\varepsilon^{2}+\frac{8}{\varepsilon}o(\sum_{k=1}^{m_{n}}|c_{nk}%
|^{2}H(|c_{nk}|^{-1})\text{ as }n\rightarrow\infty\text{.}%
\end{gather*}
Now we take into account condition (\ref{gen}) and obtain $\mathbb{E}%
(\max_{1\leq i\leq m_{n}}|c_{ni}\xi_{i}^{\prime}|^{2})\rightarrow0$ by letting
first $n\rightarrow\infty$ followed by $\varepsilon\rightarrow0.$

In order to verify the item $(b)$ of Theorem \ref{MARTCLT}, we have to study
the limit in probability of $\sum_{k=1}^{m_{n}}c_{nk}^{2}(\xi_{k}^{\prime
})^{2}.$ We start from the decomposition
\begin{gather*}
\sum_{k=1}^{m_{n}}c_{nk}^{2}(\xi_{k}^{\prime})^{2}=\sum_{k=1}^{m_{n}}%
c_{nk}^{2}\xi_{k}^{2}I(|c_{nk}\xi_{k}|\leq1)+\sum_{i=1}^{m_{n}}c_{ni}%
^{2}\mathbb{E}_{i-1}^{2}(\xi_{i}I(|c_{ni}\xi_{i}|\leq1))\\
+2\sum_{k=1}^{m_{n}}c_{nk}^{2}\xi_{k}I(|c_{nk}\xi_{k}|\leq1)\mathbb{E}%
_{k-1}(\xi_{k}I(|c_{nk}\xi_{k}|>1))=A+B+2C\text{.}%
\end{gather*}
We shall show that it is enough to analyze the first term by the following
simple argument.

We discuss now the term $B$. \ By the fact that $c_{ni}^{2}(\mathbb{E}%
_{i-1}(\xi_{i}I(|c_{ni}\xi_{i}|\leq1)))^{2}\leq1$ a.s. we obtain%
\[
B\leq\sum_{i=1}^{m_{n}}|c_{ni}\mathbb{E}_{i-1}(\xi_{i}I(|c_{ni}\xi_{i}%
|\leq1))|\ \text{\ a.s.}%
\]
By the martingale property $\mathbb{E}_{i-1}(\xi_{i}I(|c_{ni}\xi_{i}%
|\leq1))=\mathbb{E}_{i-1}(\xi_{i}I(|c_{ni}\xi_{i}|>1))$ a.s. and taking into
account (\ref{negl}) and the properties of conditional expectation,
\begin{align}
\mathbb{E}(B)  &  \leq\mathbb{E}(\sum_{i=1}^{m_{n}}|c_{ni}\mathbb{E}_{i-1}%
(\xi_{i}I(|c_{ni}\xi_{i}|>1))|)\label{negl2}\\
&  \leq\sum_{i=1}^{m_{n}}|c_{ni}|\mathbb{E(}|\xi_{0}|I(|c_{ni}\xi
_{0}|>1))=o(1)\text{ as }n\rightarrow\infty\text{.}\nonumber
\end{align}
Then, by Cauchy Schwarz inequality, condition (\ref{gen}) and (\ref{negl2}),
we get
\[
\mathbb{E}(C)\leq(\mathbb{E}(A)\mathbb{E}(B))^{1/2}\rightarrow0\text{ as
}n\rightarrow\infty\text{ .}%
\]
By these arguments, the limit in probability of $\sum_{k=1}^{m_{n}}c_{nk}%
^{2}(\xi_{k}^{\prime})^{2}$ coincides to the limit of $\sum_{k=1}^{m_{n}%
}c_{nk}^{2}\xi_{k}^{2}I(|c_{nk}\xi_{k}|\leq1).$ Notice now that
\[
\sum_{k=1}^{m_{n}}c_{nk}^{2}\xi_{k}^{2}I(|c_{nk}\xi_{k}|\leq1)-1=(\sum
_{k=1}^{m_{n}}c_{nk}^{2}\xi_{k}^{2}-1)+\sum_{k=1}^{m_{n}}c_{nk}^{2}\xi_{k}%
^{2}I(|c_{nk}\xi_{k}|>1)\text{.}%
\]
By the martingale inequality in Burkholder\ (1966) we obtain that there is a
positive constant $c$, such that for any $\varepsilon>0$
\[
\mathbb{P}(\sum_{k=1}^{m_{n}}c_{nk}^{2}(\xi_{k}^{\prime\prime})^{2}%
>\varepsilon/2)\leq\frac{c}{\varepsilon^{2}}\mathbb{E}|\sum_{k=1}^{m_{n}%
}c_{nk}(\xi_{k}^{\prime\prime})|
\]
which, combined with (\ref{negl}) and the arguments in (\ref{negl2}) gives%
\begin{equation}
\sum_{k=1}^{m_{n}}c_{nk}^{2}\xi_{k}^{2}I(|c_{nk}\xi_{k}|>1)\Rightarrow
0\text{.} \label{negl5}%
\end{equation}
We just have to take into account condition (\ref{convprob}) to conclude
\[
\sum_{k=1}^{m_{n}}c_{nk}^{2}\xi_{k}^{2}I(|c_{nk}\xi_{k}|\leq1)\Rightarrow
1\text{,}%
\]
and then, by the above considerations, $\sum_{k=1}^{m_{n}}c_{nk}^{2}(\xi
_{k}^{\prime})^{2}\Rightarrow1$.

The conclusion of this theorem follows by using the Theorem \ref{MARTCLT}.
$\diamondsuit$

\bigskip

\bigskip{\Large Proof of Proposition \ref{mix}.}

\bigskip

We assume that (M$_{1}$) holds. Because of (\ref{negl5}) \ it is enough to
prove
\[
\sum_{k=1}^{m_{n}}c_{nk}^{2}\xi_{k}^{2}I(|c_{nk}\xi_{k}|\leq1)\Rightarrow1.
\]
We start by computing the variance$~$\ of $\sum_{k=1}^{m_{n}}c_{nk}^{2}\xi
_{k}^{2}I(|c_{nk}\xi_{k}|\leq1)$ and, taking into account the variables are
identically distributed, we majorate the covariances by using the coefficients
$\psi_{i}$ defined in Condition (M$_{1}$).
\begin{gather*}
var(\sum_{k=1}^{m_{n}}c_{nk}^{2}\xi_{k}^{2}I(|c_{nk}\xi_{k}|\leq1))\leq
\sum_{k=1}^{m_{n}}c_{nk}^{4}\mathbb{E(}\xi_{0}^{4}I(|c_{nk}\xi_{0}|\leq1))\\
+2\sum_{i=1}^{m_{n}-1}\psi_{i}\sum_{k=1}^{m_{n}-i}c_{nk}^{2}c_{n,k+i}%
^{2}\mathbb{E(}\xi_{0}^{2}I(|c_{nk}\xi_{0}|\leq1))\mathbb{E(}\xi_{0}%
^{2}I(|c_{n,k+i}\xi_{0}|\leq1))=D+2E\text{.}%
\end{gather*}
By the item 4 of Lemma \ref{csorgo},%
\[
D=\sum_{k=1}^{m_{n}}c_{nk}^{4}\mathbb{E(}\xi_{0}^{4}I(|c_{nk}\xi_{0}%
|\leq1))=o(\sum_{k=1}^{m_{n}}c_{nk}^{2}H(|c_{nk}|^{-1}))=o(1)\text{ as
}n\rightarrow\infty\text{.}%
\]
In order to estimate the second term we split the sum in two, one up to $h$
and another after $h$, where $h$ is an integer.
\begin{align*}
E  &  =\sum_{i=1}^{m_{n}-1}\psi_{i}\sum_{k=1}^{m_{n}-i}c_{nk}^{2}c_{n,k+i}%
^{2}H(|c_{nk}|^{-1})H(|c_{n,k+i}|^{-1})\leq\\
&  h\max_{1\leq i\leq h}\sum_{k=1}^{m_{n}-i}c_{nk}^{2}c_{n,k+i}^{2}%
H(|c_{nk}|^{-1})H(|c_{nk+i}|^{-1})\\
&  +\max_{h\leq i\leq m_{n}}\psi_{i}\sum_{i=1}^{m_{n}-1}\sum_{k=1}^{m_{n}%
-i}c_{nk}^{2}c_{n,k+i}^{2}H(|c_{nk}|^{-1})H(|c_{n,k+i}|^{-1})\\
&  \leq h\max_{1\leq k\leq m_{n}}c_{nk}^{2}H(|c_{nk}|^{-1})+\max_{h\leq i\leq
m_{n}}\psi_{i}\text{.}%
\end{align*}
By letting $n\rightarrow\infty$ and then $h\rightarrow\infty$, by taking into
account conditions (\ref{gen}) and \textbf{(}M$_{1}$)\textbf{, }we obtain
$E\rightarrow0$ as $n\rightarrow\infty$. Then, clearly by the above
considerations $var(\sum_{k=1}^{m_{n}}c_{nk}^{2}(\xi_{k}^{\prime}%
)^{2})\rightarrow0$ that further implies by condition (\ref{gen}) that
(\ref{convprob}) holds under (M$_{1}$).

The proof of this proposition under (M$_{2}$) is similar. Because $c_{ni}%
^{2}\mathbb{E(}\xi_{i}^{2}I(|c_{ni}\xi_{i}|\leq1))\leq1$,
\begin{gather*}
var(\sum_{k=1}^{m_{n}}c_{nk}^{2}\xi_{k}^{2}I(|c_{nk_{k}}\xi|\leq1)\leq
\sum_{k=1}^{m_{n}}c_{nk}^{4}\mathbb{E(}\xi_{0}^{4}I(|c_{nk}\xi_{0}|\leq1))\\
+2\sum_{i=1}^{m_{n}-1}\sum_{k=i+1}^{m_{n}}\varphi_{k-i}c_{nk}^{2}%
\mathbb{E(}\xi_{0}^{2}I(|c_{nk}\xi_{0}|\leq1))\text{ .}%
\end{gather*}
The first term in the right hand side is treated as before. For $h$ fixed we
easily obtain
\begin{align*}
\sum_{i=1}^{m_{n}-1}\sum_{k=i+1}^{m_{n}}\varphi_{k-i}c_{nk}^{2}\mathbb{E}%
\xi_{0}^{2}I(|c_{nk}\xi_{0}|  &  \leq1)\leq h\max_{1\leq k\leq m_{n}}%
c_{nk}^{2}H(|c_{nk}|^{-1})\sum_{k=1}^{m_{n}}\varphi_{k}\\
+\sum_{i=h}^{\infty}\varphi_{i}\sum_{k=1}^{m_{n}}c_{nk}^{2}\mathbb{E}\xi
_{0}^{2}I(|c_{nk}\xi_{0}|  &  \leq1)\text{,}%
\end{align*}
and the result follows by (\ref{gen}) and condition (M$_{2}$) by letting first
$n\rightarrow\infty$ followed by $h\rightarrow\infty.$ $\diamondsuit$

\bigskip

{\Large Proof of Proposition \ref{Exist}}

\bigskip

By the three series theorem for martingales, (Theorem 2.16 in Hall and Heyde,
1980) $X_{0}$ exists in the almost sure sense if and only if:

\begin{enumerate}
\item $\sum_{i\ }^{\ }\mathbb{P}(|a_{i}\mathbb{\xi}_{i}|>1|\mathcal{F}%
_{i-1})<\infty$ a.s.,

\item $\sum_{i\ }^{\ }\mathbb{E}(a_{i}\mathbb{\xi}_{i}I(|a_{i}\mathbb{\xi}%
_{i}|\leq1)|\mathcal{F}_{i-1})$ converges a.s.,

\item $\sum_{i\ }^{\ }Var(a_{i}\mathbb{\xi}_{i}I(|a_{i}\mathbb{\xi}_{i}%
|\leq1|\mathcal{F}_{i-1}))<\infty$ \ a.s.
\end{enumerate}

Notice that, by taking into account Convention 1, the fact that the variables
are identically distributed and item 2 in Lemma \ref{csorgo} from the
appendix,
\[
\sum_{i}\mathbb{P}(|a_{i}\mathbb{\xi}_{i}|>1)=\sum_{i}\mathbb{P}%
(|a_{i}\mathbb{\xi}_{0}|>1)=\sum_{i}a_{i}^{2}o(H(|a_{i}|^{-1}))<\infty\text{ }%
\]
and this easily implies 1.

Then by item 3 of Lemma \ref{csorgo} and again by the fact that the variables
are identically distributed,
\[
|\sum_{i}\mathbb{E(}a_{i}\mathbb{\xi}_{i}I(|a_{i}\mathbb{\xi}_{i}|\leq
1)|\leq\sum_{i}|a_{i}|\mathbb{E(}|\mathbb{\xi}_{0}|I(|a_{i}\mathbb{\xi}%
_{0}|>1))=\sum_{i}a_{i}^{2}o(H(|a_{i}|^{-1})<\infty\text{.}%
\]
This implies%
\begin{equation}
\sum_{i}\mathbb{E(}|a_{i}\mathbb{\xi}_{i}|I(|a_{i}\mathbb{\xi}_{i}%
|\leq1)|\mathcal{F}_{i-1})<\infty\ \text{\ a.s.} \label{negl3}%
\end{equation}
and 2. follows.

Finally,%
\[
\sum_{i}\mathbb{E(}a_{i}^{2}\mathbb{\xi}_{i}^{2}I(|a_{i}\mathbb{\xi}_{i}%
|\leq1))=\sum_{i}a_{i}^{2}\mathbb{E(\xi}_{0}^{2}I(|a_{i}\mathbb{\xi}_{0}%
|\leq1))=\sum_{i}a_{i}^{2}H(|a_{i}|^{-1}).
\]
Then%
\[
\sum_{i}a_{i}^{2}\mathbb{E(\xi}_{i}^{2}I(|a_{i}\mathbb{\xi}_{i}|\leq
1)|\mathcal{F}_{i-1})<\infty\text{ a.s.}%
\]
and together with (\ref{negl3}) gives 3. For the i.i.d. case the proof is
similar and it is based on the i.i.d. version of the three series theorem.
$\diamondsuit$

\bigskip

{\Large Proof of Theorem \ref{linCLT}}

\bigskip

We start by rewriting $S_{n}$ as in relation (\ref{defS}) by changing the
order of summation.
\[
S_{n}=\sum_{k=1}^{n}X_{k}=\sum_{j=-\infty}^{\infty}(\sum_{k=1}^{n}a_{k+j}%
)\xi_{j}=\sum_{j=-\infty}^{\infty}b_{nj}\xi_{j}\text{.}%
\]
We shall verify the conditions of Theorem \ref{genthm}.

According to Corollary (\ref{CLTnorm}) it is sufficient to show that
\[
\sup_{k}\frac{b_{nk}^{2}}{D_{n}^{2}}H(\frac{D_{n}}{|b_{nk}|})\rightarrow
0\;\;\text{as}\;n\rightarrow\infty\text{.}%
\]
where
\[
D_{n}=\inf\{x\geq1:\sum_{k}\frac{b_{nk}^{2}}{x^{2}}H(\frac{x}{|b_{nk}|}%
)\leq1\}\text{.}%
\]
Notice that condition (\ref{coeff0}) implies $\sum_{k}a_{k}^{2}<\infty.$
Therefore we can apply the argument from Peligrad and Utev (1997, pages
448-449) and obtain
\begin{equation}
\sup_{k}\frac{b_{nk}^{2}}{\sum_{i}b_{ni}^{2}}\rightarrow0\text{ }\ \text{as
}n\rightarrow\infty\text{,} \label{negl4}%
\end{equation}
since we imposed $\sum_{k}b_{nk}^{2}\rightarrow\infty.$ Notice that by taking
into account Convention 2 we obviously have
\[
\frac{1}{\sum_{k}b_{nk}^{2}}\sum_{k}b_{nk}^{2}H\left(  \frac{(\sum_{k=1}%
b_{nk}^{2})^{1/2}}{|b_{nk}|}\right)  \geq1\text{.}%
\]
By the definition of $D_{n}$, this implies
\[
D_{n}^{2}\geq\sum_{k}b_{nk}^{2}\text{,}%
\]
whence, by (\ref{negl4})
\[
\sup_{k}\frac{b_{nk}^{2}}{D_{n}^{2}}\rightarrow0\text{ }\ \text{as
}n\rightarrow\infty\text{.}%
\]
Now, by the properties of slowly varying functions, for $\varepsilon>0$, we
know that $H(x)=o(x^{\varepsilon})$ as $x\rightarrow\infty.$ We then obtain
\[
\sup_{k}\frac{b_{nk}^{2}}{D_{n}^{2}}H(\frac{D_{n}}{|b_{nk}|})=o(1)\text{ as
}n\rightarrow\infty\text{.}%
\]
This completes the proof of this Theorem. $\diamondsuit$

\bigskip

{\Large Proof of Theorem \ref{MIXING}}

\bigskip

In order to prove Theorem \ref{MIXING} we start by the truncation argument \
\[
\xi_{ni}^{\prime}=\xi_{i}I(|c_{ni}\xi_{i}|\leq1)-\mathbb{E}(\xi_{i}%
I(|c_{ni}\xi_{i}|\leq1))
\]
and%
\[
\xi_{ni}^{\prime\prime}=\xi_{i}I(|c_{ni}\xi_{i}|>1)-\mathbb{E}(\xi
_{i}I(|c_{ni}\xi_{i}|>1))\text{.}%
\]
\textbf{ }First, by relation (\ref{negl}), we have
\begin{equation}
E|\sum_{k=1}^{m_{n}}c_{nk}\xi_{k}^{\prime\prime}|\rightarrow0\text{,}
\label{negl6}%
\end{equation}
whence, by Theorem 3.1 in Billingsley (1999), the proof is reduced to studying
the asymptotic behavior of $\sum_{k=1}^{m_{n}}c_{nk}\xi_{k}^{\prime}$.

According to Theorem 4.1 and Theorem 5.5 in Utev (1990), given for convenience
in the appendix (Theorem \ref{Utev}), we have only to verify the Lindeberg's
condition. Denote $(\sigma_{n}^{\prime})^{2}=\mathbb{E}(\sum_{k=1}^{m_{n}%
}c_{nk}\xi_{k}^{\prime})^{2}.$ Our conditions on the mixing coefficients
allows us by relation (\ref{VA}) to bound $(\sigma_{n}^{\prime})^{2}$ above
and below by the sum of squares of the variance of individual summands; so we
can find two positive constants $C_{1}<C_{2}$ such that%
\begin{equation}
C_{1}\sum_{i=1}^{m_{n}}c_{ni}^{2}H(|c_{ni}|^{-1})\leq(\sigma_{n}^{\prime}%
)^{2}\leq C_{2}\sum_{i=1}^{m_{n}}c_{ni}^{2}H(|c_{ni}|^{-1}). \label{var}%
\end{equation}
Lindeberg's condition is satisfied because Lyapunov's condition is satisfied.
To see this, by (\ref{var}) we have
\[
\frac{1}{(\sigma_{n}^{\prime})^{2}}\sum_{i=1}^{m_{n}}c_{ni}^{3}\mathbb{E(}%
|\xi_{0}|^{3}I(|c_{ni}\xi|\leq1))\leq\frac{1}{(\sigma_{n}^{\prime})^{2}}%
\sum_{i=1}^{m_{n}}c_{ni}^{2}o(H(|c_{ni}|^{-1}))=o(1)\text{.}%
\]
The central limit theorem follows and we have
\begin{equation}
(\sum_{k=1}^{m_{n}}c_{nk}\xi_{k}^{\prime})/\sigma_{n}^{\prime}\Rightarrow
N(0,1)\;\;\text{as}\;\;n\rightarrow\infty\text{.} \label{CLTprime}%
\end{equation}
Now, $(\sum_{k=1}^{m_{n}}c_{nk}\xi_{k}^{\prime})/\sigma_{n}^{\prime}$ is
uniformly integrable and by Theorem 3.4 in Billingsley (1999), this implies
that $\ \mathbb{E(}|\sum_{k=1}^{m_{n}}c_{nk}\xi_{k}^{\prime}|/\sigma
_{n}^{\prime})\rightarrow\sqrt{2/\pi}$. By taking now into account
(\ref{negl6}) we obtain $\ \mathbb{E(}|\sum_{k=1}^{m_{n}}c_{nk}\xi_{k}%
|/\sigma_{n}^{\prime})\rightarrow\sqrt{2/\pi}$, that combined with
(\ref{CLTprime}) and (\ref{negl6}), leads to the conclusion of this theorem.
$\diamondsuit$

\section{Appendix}

We mention first a lemma which contains some equivalent formulation for
variables in the domains of attraction of normal law. It is Lemma 1 in
Cs\"{o}rg\H{o} et al (2003).

\begin{lemma}
\label{csorgo} Let $H(x):=\mathbb{E(}X^{2}I(|X|\leq x))$. The following
statements are equivalent.
\end{lemma}

\begin{enumerate}
\item $H(x)$ is a slowly varying function at $\infty$;

\item $\mathbb{P}(|X|>x)=o(x^{-2}H(x))$;

\item $\mathbb{E(}|X|I(|X|>x))=o(x^{-1}H(x))$;

\item $\mathbb{E(}|X|^{\alpha}I(|X|\leq x))=o(x^{\alpha-2}H(x))$ for
$\alpha>2$.
\end{enumerate}

The following theorem is a simplified version of Theorem 3.2 from Hall and
Heyde (1980). See also Relation 2.18 and (M) in Gaenssler and Haeusler (1986).

\begin{theorem}
\label{MARTCLT}Let $(D_{ni})_{1\leq i\leq k_{n}}$ be a square integrable
martingale difference and $(\mathcal{F}_{i})~$a filtration of sigma algebras
such that for each, $n$ and $1\leq i\leq k_{n}$, $D_{ni}$ is $\mathcal{F}_{i}$
measurable. Suppose that\newline(a) $\max_{1\leq i\leq k_{n}}|D_{ni}|$ is
uniformly integrable;\newline(b) $\sum_{i=1}^{k_{n}}D_{ni}^{2}\Rightarrow1$ as
$n\rightarrow\infty;$ \newline Then $S_{n}\;\Rightarrow\;N(0,1)$ as
$n\rightarrow\infty$ where $S_{n}=\sum_{i=1}^{k_{n}}X$\emph{${_{ni}}$.}
\end{theorem}

By combining Theorem 4.1 and Theorem 5.5 in Utev (1990) we formulate a general
result for triangular arrays of $\rho$-mixing random variables.

\begin{theorem}
\label{Utev}Assume $(Y_{nk})_{1\leq k\leq m_{n}}$\emph{ }$\ $is a triangular
array of centered random variables with finite second moment such that
$\rho(1)<1$ and $\sum_{k}\rho(2^{k})<\infty$ where
\[
\rho(k)=\sup_{s\geq1,n\geq1}\rho(\sigma(Y_{ni};i\leq s),\sigma(Y_{ni};i\geq
s+k).
\]
Set $S_{n}=\sum_{i=1}^{m_{n}}Y_{ni}$\emph{ }and denote $var(S_{n})=\sigma
_{n}^{2}.$ Then, there are two positive constants $C_{1}$ and $C_{2}$ such
that%
\begin{equation}
C_{1}\sum_{i=1}^{m_{n}}\mathbb{E}(Y_{ni}^{2})\leq\sigma_{n}^{2}\leq C_{2}%
\sum_{i=1}^{m_{n}}\mathbb{E}(Y_{ni}^{2}). \label{VA}%
\end{equation}
Moreover if Lindeberg's condition is satisfied:%
\[
\frac{1}{\sigma_{n}^{2}}\sum_{i=1}^{m_{n}}\mathbb{E}(Y_{ni}^{2}I(|Y_{ni}%
|>\varepsilon\sigma_{n}))\rightarrow0\text{,}%
\]
then
\[
\frac{1}{\sigma_{n}}\sum_{i=1}^{m_{n}}Y_{ni}\Rightarrow N(0,1).
\]

\end{theorem}

\section{Acknowledgement}

The authors are grateful to the referees for their useful suggestions which
improved the presentation of this paper.


\begin{thebibliography}{99}                                                                                               %


\bibitem {ag}Araujo, A. and Gin\'{e}, E. (1980). \textit{The Central Limit
Theorem for Real and Banach Valued Random Variables}. Wiley Series in
Probability and Mathematical Statistics. John Wiley \& Sons, New York-Chichester-Brisbane.

\bibitem {ber}Beran, J. (1994). \textit{Statistics for long-memory processes}.
Monographs on Statistics and Applied Probability, \textbf{61}. Chapman and
Hall, New York.

\bibitem {Bill}Billingsley, P. (1999). \textit{Convergence of Probability
measures}. Second edition, Wiley, New York.

\bibitem {BR}Bradley, R. C. (1988). A Central Limit Theorem for Stationary
$\rho$-Mixing Sequences with Infinite Variance. \textit{Ann. Probab}.
\textbf{16} 313-332.

\bibitem {Br}Bradley, R. C. (2007). \textit{Introduction to strong mixing
conditions}. Volumes 1-3, Kendrick Press.

\bibitem {Bu}Burkholder, D. L. (1966). Martingale transforms. \textit{Ann.
Math. Statist}. \textbf{37} 1494-1504.

\bibitem {Cs}Cs\"{o}rg\H{o}, M., Szyszkowicz, B. and Wang, Q. (2003).
Donsker's theorem for self-normalized partial sums processes, \textit{Ann.
Probab}. \textbf{31} 1228-1240.

\bibitem {Feller}Feller, W. (1966). \textit{An Introduction to Probability
Theory and Its Applications} Vol. 2, Willey, New York.

\bibitem {GH}Gaenssler, P. and Haeusler, E. (1986). On martingale central
limit theory. \textit{Dependence in Probability and Statistics},
Birkh\"{a}user, Boston, 303-334.

\bibitem {Gine}Gin\'{e}, E., G\"{o}tze, F. and Mason, D. M. (1997). When is
the Student t--statistic asymptotically standard normal? \textit{Ann. Probab.}
\textbf{25}, 1514--1531.

\bibitem {Go}Gordin, M. I. (1969). The central limit theorem for stationary
processes, \textit{Soviet. Math. Dokl.}\emph{ }\textbf{10}, 1174--1176.

\bibitem {HH}Hall P. and Heyde, C. C. (1980). \textit{Martingale Limit Theory
and Its Applications}, Academic Press.

\bibitem {J}Jara, M., Komorowski, T. and Olla, S. (2009). Limit theorems for
additive functionals of a Markov Chain. \textit{Ann. Appl. Probab.} 2009,
\textbf{19}, 2270--2300

\bibitem {}Knight, K. (1991). Limit theory for M-estimates in an integrated
infinite variance process. \textit{Econometric Theory} \textbf{7,} 200-212.

\bibitem {Kulik}Kulik, R. (2006). Limit theorems for self-normalized linear
processes. \textit{Statistics and Probability Letters} \textbf{76,} 1947-1953.

\bibitem {IL}Ibragimov, I. A. and Linnik, Yu. V. (1971). \textit{Independent
and Stationary Sequences of Random Variables}, Wolters, Groningen.

\bibitem {Mason}Mason, D.M. (2005). The asymptotic distribution of
self-normalized triangular arrays. \textit{J. Theor. Probab.} \textbf{18}, 853-870.

\bibitem {m}Mikosch, T. Gadrich, T. Kliippelberg C. and Adler, R. J. (1995).
Parameter estimation for ARMA models with infinite variance innovations.
\textit{Annals of Statistics} \textbf{23,} 305-326.

\bibitem {PU1}Peligrad, M. and Utev, S. (1997). Central limit theorem for
linear processes. \textit{Ann. Probab}. \textbf{25}, 443-456.

\bibitem {PU2}Peligrad, M. and Utev, S. (2006). Central limit theorem for
stationary linear processes. \textit{Ann. Probab}. \textbf{34}, 1608-1622.

\bibitem {PS}Peligrad, M. and Sang, H. (2010). Asymptotic Properties of
Self-Normalized Linear Processes with Long Memory. arXiv:1006.1572

\bibitem {Ro}Robinson, P. M. (1997). Large-sample inference for non parametric
regression with dependent errors. \textit{Ann. Statist. }\textbf{25}, 2054-2083.

\bibitem {Shao}Shao, Q. (1993). On the invariance principle for$\rho$-mixing
sequence of random variables with infinite variance. \textit{Chinese Ann.
Math.} Ser. B. \textbf{14}, 27-42.

\bibitem {Tir}Tyran-Kami\'{n}ska, M. (2010). Functional limit theorems for
linear processes in the domain of attraction of stable laws .
\textit{Statistics and Probability Letters} \textbf{77}, 1535-1541.

\bibitem {Ut}Utev, S. A. (1990). Central limit theorem for dependent random
variables. \textit{Prob. Theory and Math. Stat.} Vol. \textbf{2}, B.
Grigelionis et al. (eds.), VSP/Mokslas. 519-528.

\bibitem {V}Voln\'{y}, D (1993). Approximating martingales and central limit
theorem for strictly stationary processe. \textit{Stoch. Proc. Appl.}
\textbf{44,} 41-74.

\bibitem {W}Wu, W. B. (2003). Additive functionals of infinite-variance moving
averages. \textit{Statistica Sinica} \textbf{13}, 1259-1267.
\end{thebibliography}
\end{document}